\documentclass[11pt,twoside,a4paper]{article}

\usepackage{graphics,epsfig}
\usepackage[margin=1in]{geometry}

\title{\LARGE \bf
Dilated Matrix Inequalities for Control Design in Systems with
Actuator Constraint\thanks{This work was supported by NSF Grant CMS-0510874}}
\date{}

\author{Solmaz Sajjadi-Kia and Faryar Jabbari
\thanks{The authors are with the Department of Mechanical and Aerospace Engineering, University of California, Irvine,
        Irvine, CA 92697, USA.
        {\tt\small ssajjadi,fjabbari@uci.edu}}%
}

\begin{document}

\maketitle \thispagestyle{empty} \pagestyle{empty}
\newtheorem{Lem}{Lemma}
\newtheorem{Thm}{Theorem}
\newtheorem{REMARK}{Remark}
\newtheorem{proof}{Proof}

\begin{abstract}

In this paper, we present a new variation of dilated matrix
inequalities (MIs) for Bounded Real MI, invariant set MI and
constraint MI, for both state and output feedback synthesis
problems. In these dilated MIs, system matrices are separated from
Lyapunov matrices to allow the use of different Lyapunov matrices
in multi-objective and robust problems. To demonstrate the benefit
of these new dilated MIs over conventional ones, they are used in
solving controller synthesis problem for systems with bounded
actuator in disturbance attenuation. It is shown that for the
resulting multi-objective saturation problem, the new form of
dilated MIs achieves an upper bound for $L_2$ gain that is less
than or equal to the upper bound estimate achieved by conventional
method.

\end{abstract}

\section{INTRODUCTION}

Most of the (Linear) Matrix Inequality ((L)MI) characterizations
in control techniques such as ${\small H_\infty}$, $H_2$, use a
quadratic Lyapunov function (${\small V=x^TPx\quad |P>0}$ ) to
develop their MIs (e.g., see ref \cite{c1}). The resulting MIs end
up having entries with the products of the Lyapunov variables and
system matrices. This causes some degree of conservatism in
multi-objective and robust problems by forcing common Lyapunov
matrices for all objectives. For example, see Ref.
\cite{c3}-\cite{c5} which make use of common Lyapunov variable in
their multi-objective problems.

\smallskip
Recently, researchers have been using matrix dilation results to
reduce this conservatism. The earlier, and the best known, of the
results obtained by these new techniques are in discrete time
settings (\cite{c7},\cite{c8}). Although a lot of effort has also
been made for the continuous-time case, it is still an open
problem, mostly due to the fact that dilations to reduce
conservatism destroy the convexity in some important cases. Some
of the very nice and
convex results in continuous time are achieved by Ebihara et al. in \cite{c9}-\cite{c11}, in case where 
the problem can be case as
\begin{equation}\label{eq::ebihara}{\small
AX+XA^T+\delta_1X+\delta_2AXA^T+X\Delta^T\Delta X<0}
\end{equation}
for a suitable choice of ${\small X}$, ${\small \Delta}$, etc. By
assigning different matrices to ${\small A}$, ${\small \delta_1}$,
${\small\delta_2}$ and ${\small \Delta}$, this general form covers
some continuous-time control problems such as stability, $\small
H_2$ and D-stability. Furthermore, references \cite{c16} and
\cite{c17} present a technique that, using a projection (or
Elimination lemma) based approach, leads to set of convex search
for several important problems; ${\small H_2}$, stability,
D-stability, etc.

\smallskip
Unfortunately, neither of the two approaches above
deal with synthesis inequalities faced in the $L_2$-gain (i.e.,
bounded real problems) or several invariant set determination
problems. Roughly, these problems result in a constant term in
(\ref{eq::ebihara}). There have been dilated MIs for such cases,
as well (\cite{c12}-\cite{c14}, among others). So far, there seems
to be two weaknesses associated with these set of results. Often,
the synthesis results are for the full state case, by exploiting a
structure that holds only in full state problem. Furthermore, they
all seem to need an additional scalar variable which enters the MI
in a  non-convex  fashion (\cite{c12}-\cite{c14}). While it seems
that such a non-convexity is inevitable, it can be addressed with
a line search since the culprit is a scalar variable.

\smallskip
Here, we present a new dilated MI, that can be applied to the
bounded real LMI, as well as to the matrix inequalities that are
used in the invariant set for peak bounded disturbance (\cite{c1})
and the constraint LMI (\cite{c1}). These new MIs are obtained
explicitly through a constructive methods, to avoid the
ambiguities that can -- at time -- accompany results based on the
Projection Lemma (see \cite{c9} for a discussion on this issue).
We also rely on a scalar variable that renders the problem
non-convex and use line searches to obtain the final result.
Fortunately, we show that the proposed approach is rather easily
extended to the output feedback problem, assuming the controller
is of the same order as the plant.

\smallskip
Finally, we study the effect of this new matrix dilation
technique in reducing conservatism in controller design for linear
systems with bounded actuators which avoid saturation. As
mentioned in \cite{c3}, this is often set as a multi-objective
problem, and can suffer the conservatism forced by common Lyapunov
matrix. We show that, these new MIs lead us to a problem whose
performance is at least equal to the one with standard MIs.

\smallskip
The system we study has the standard model
\begin{equation}\label{eq::OL_system}{\small
\left\{ \begin{array}{l}
\dot x_p  = Ax_p+B_1w+B_2u\\
z =C_1x_p+D_{11}w+D_{12}u\\
y=C_2x_p+D_{21}w
\end{array} \right.}
\end{equation}
with the closed loop
\begin{equation}\label{eq::CL_system}{\small
\left\{ \begin{array}{ll}
\dot x = A_{cl}x+B_{cl}w\\
z=C_{cl}x+D_{cl}w\\
\end{array} \right.}
\end{equation}
where the details differ in the state feedback and output feedback
cases. The Transfer function of this system is ${\small
T_{zw}(s)=C_{cl}(sI-A_{cl})^{-1}B_{cl}+D_{cl}}$. We use ${\small
\mathbf{He}(A)}$ as short notation for ${\small A+A^T}$ to save
space. The rest of the notations throughout the paper follow
standard practices.

\section{Dilated MIs for Some Practical Design Specifications}
In this section, we derive equivalent dilated MIs for standard
Bounded real LMI, invariant set MI and constraint LMI.  In the new
MIs, the system matrices and the Lyapunov variable are decoupled.

\subsection{A Dilated LMI for $L_2$ Gain}

\begin{Lem}[Bounded Real Lemma \cite{c1}]\label{Lem_Bounded_Real}
$A_{cl}$ is stable with ${\small
\|T_{zw}(s)\|_\infty<\gamma_{con}}$ if and only if ${\small
\bar{\sigma}(D_{cl})<\gamma_{con}}$ and there exists $Q_1>0$ such
that
\begin{equation}\label{eq::Q_Perf}{\small
\left( \begin{array}{ccc}
A_{cl}Q_1+Q_1A_{cl}^T & B_{cl} & Q_1C_{cl}^T\\
\star & -\gamma_{con} I & D_{cl}^T\\
\star & \star & -\gamma_{con} I
\end{array} \right)<0}
\end{equation}
\end{Lem}

\bigskip
\begin{Thm}\label{TH_DMI_L2}
The closed-loop system (\ref{eq::CL_system}) is stable and its
$L_2$ gain is less than $\gamma_{new}$ if there exist a positive
constant $0<\epsilon_1<1$, and square matrices $X_1>0$ and $G_1$
which satisfy: \setlength\arraycolsep{0.7pt}
\begin{equation}\label{eq::eq_perf}{\small
\left( \begin{array}{cccc}
X_1 & B_{cl} & 0 & -X_1\\
\star & -\gamma_{new} I & D_{cl}^T & 0 \\
\star & \star & -\gamma_{new} I & 0\\
\star & \star & \star & 0 \end{array}
\right)+\mathbf{He}(\mathcal{Q}^TG_1 \mathcal{P})<0}
\end{equation}
\end{Thm}\medskip
where $${\small \mathcal{P}=[I\quad 0 \quad 0 \quad -2\epsilon_1
I];\quad \mathcal{Q}=[(A_{cl}^T-\frac{1}{2}I) \quad 0 \quad
C_{cl}^T\quad I].}$$

\begin{proof}
Proof is through showing MIs (\ref{eq::Q_Perf}) and
(\ref{eq::eq_perf}) are equivalent. Suppose that MI
(\ref{eq::eq_perf}) holds. Consider the explicit bases of
nullspaces of $\mathcal{P}$ and $\mathcal{Q}$
\begin{displaymath}{\small \setlength\arraycolsep{4pt}
\mathcal{N}_{\mathcal{P}}=\left( \begin{array}{ccc}
I & 0 & 0\\
0 & I & 0 \\
0 & 0 & I\\
\frac{1}{2\epsilon_1}I & 0 & 0\end{array} \right),\,
\mathcal{N}_{\mathcal{Q}}=\left(
\begin{array}{ccc}
I & 0 & 0\\
0 & I & 0 \\
0 & 0 & I\\
-A_{cl}^T+\frac{1}{2}I & 0 & -C_{cl}^T\end{array} \right)\,\,}
\end{displaymath}
By multiplying $\mathcal{N}_{\mathcal{Q}}$ and its transpose from
the right and left sides respectively, and considering
$\mathcal{Q}\mathcal{N}_{\mathcal{Q}}=0$, inequality
(\ref{eq::eq_perf}) becomes
\begin{equation}\label{eq::Th1_1}{\small
\left( \begin{array}{ccc}
A_{cl}X_1+X_1A^T_{cl} & B_{cl} & X_1C_{cl}^T\\
\star & -\gamma_{new} I & D_{cl}^T\\
\star & \star & -\gamma_{new} I
\end{array} \right)<0}
\end{equation}
Now if we define $X_1=Q_1$ and $\gamma_{new}=\gamma_{con}$, it is
clear that (\ref{eq::Q_Perf}) holds.
Next, multiplying (\ref{eq::eq_perf}) from right and left by
$\mathcal{N}_{\mathcal{P}}$ and its transpose respectively gives

\begin{equation}\label{eq::Th1_3}{\small
\left( \begin{array}{ccc}
(1-\frac{1}{\epsilon_1})X_1 & B_{cl} & 0\\
\star & -\gamma_{new} I & D_{cl}^T\\
\star & \star & -\gamma_{new} I
\end{array} \right)<0}
\end{equation}
This inequality implies that (\ref{eq::eq_perf}) can have a
solution only for $0<\epsilon_1<1$.

\medskip On the other hand, suppose that (\ref{eq::Q_Perf}) holds
with $Q_1>0$. Note that (\ref{eq::Q_Perf}) can be rewritten as
(\ref{eq::Th1_1}) by defining $Q_1=X_1$ and
$\gamma_{con}=\gamma_{new}$.
Since $X_1>0$, for any $\bar{\epsilon}_1>0$ we have
\begin{equation}\label{eq::Th1_2}{\small
\mathcal{R}^T(4\bar{\epsilon_1} X_1)^{-1}\mathcal{R}\geq 0}
\end{equation}
where ${\small \mathcal{R}=[-2\bar{\epsilon_1}
X_1(A_{cl}^T-\frac{I}{2})\quad 0\quad -2\bar{\epsilon_1} X_1
C_{cl}^T]}$. Since the right hand side of the inequality is of
order $\bar{\epsilon_1}$, it is possible to find a sufficiently
small $\bar{\epsilon_1}>0$ which for this $\bar{\epsilon_1}$, and
any $\epsilon_1<\bar{\epsilon_1}$, the following holds
$${\small [\textrm{left side of (\ref{eq::Th1_1})}]+[\textrm{left side of
(\ref{eq::Th1_2})}]<0.}$$

\noindent Applying the Schur complement to this inequality leads
to \setlength\arraycolsep{4pt}
\begin{displaymath}{\small
\left( \begin{array}{cccc}
\Pi & B_{cl} & X_1C_{cl}^T & -2\epsilon_1 (A_{cl}-\frac{1}{2}I)X_1\\
\star & -\gamma_{new} I & D_{cl} & 0 \\
\star & \star & -\gamma_{new} I & -2\epsilon_1 C_{cl}X_1 \\
\star & \star & \star & -4\epsilon_1 X_1
\end{array} \right)<0}
\end{displaymath}
where ${\small \Pi=\mathbf{He}((A_{cl}-\frac{1}{2}I)X_1)+X_1}$. By
choosing ${\small G_1=G_1^T=X_1}$, this inequality can be written
as (\ref{eq::eq_perf}); i.e., satisfaction of (\ref{eq::Q_Perf})
leads to a specific choice for the matrices $G_1$ and $X_1$ that
satisfy (\ref{eq::eq_perf}), with the same performance estimate.
\end{proof}

\medskip
\subsection{A Dilated MI for the Invariant set}
\begin{Lem}[
Invariant Set MI \cite{c1}] $\mathcal{E}=\{x|x^TPx<w^2_{max}\}$ is
a reachable set (invariant set) for the LTI system
(\ref{eq::CL_system}) exposed to a peak bounded disturbance
$w(t)^Tw(t)\leq w^2_{max}$  if there exist a scalar
$\alpha_{con}>0$ and $\small Q_2=P^{-1}>0$ such that the following
MI is feasible
\begin{equation}\label{eq::Q_Invariant}{\small
\left( \begin{array}{cc}
A_{cl}Q_2+Q_2A^T_{cl}+\alpha_{con} Q& B_{cl}\\
\star & -\alpha_{con} I
\end{array} \right)<0}
\end{equation}
\end{Lem}

\bigskip
\begin{Thm}\label{TH_DMI_Inva}
Inequality (\ref{eq::Q_Invariant}) is feasible if and only if
there exist a square matrix $G_2$, a constant $0 < \epsilon_2 < 1$
and $X_2 = P^{-1}> 0$ which satisfy the following condition for
some $\alpha_{new} > 0$:
\begin{equation}\label{eq::eq_Invariant}{\small
\left( \begin{array}{cccc}
X_2 & B_{cl} & -X_2\\
\star & -\alpha_{new} I & 0 \\
\star & \star & 0 \end{array} \right)\,+\,
\mathbf{He}(\mathcal{Q}^TG_2 \mathcal{P})<0}
\end{equation}
where ${\small \mathcal{P}=[I\quad 0 \quad -2\epsilon_2 I]}$ and
${\small \mathcal{Q}=[(\hat{A}^T-\frac{1}{2}I) \quad 0 \quad I]}$,
with ${\small\hat{A}=A_{cl}+\frac{\alpha_{new}}{2}I}$.
\end{Thm}
\medskip
\begin{proof}
Proof is similar to the Theorem 1. In this case $\mathcal{R}$ in
inequality (\ref{eq::Th1_2}) is
$${\small \mathcal{R} = [-2\bar{\epsilon_2}
X_2(\hat{A}^T-\frac{1}{2}I)\quad 0]} $$ \end{proof}

\subsection{A Dilated MI for Constraint LMI}
Constraint LMI (\cite{c1}): Suppose that $y=Kx$ and $x \in
\{x|x^TPx<w_{max}^2\}$. Then, $\|y\|^2$ will be less than or equal
to $u_{lim}^2$ if the following LMI holds for $Q_3=P^{-1}$
\begin{equation}\label{eq::Q_constraint}{\small
\left( \begin{array}{cc}
-Q_3& -Q_3K^T\\
\star & -\frac{u_{lim}^2}{\omega_{max}^2} I
\end{array} \right)<0}
\end{equation}

\begin{Thm}\label{TH_DMI_cons}
The matrix inequality (\ref{eq::Q_constraint}) is feasible if and
only if the following inequality is feasible for some
$0<\epsilon_3<1$
and $X_3>0$: 
\setlength\arraycolsep{4pt}
\begin{equation}\label{eq::eq_constraint}{\small
\left( \begin{array}{ccc}
X_3 & 0 & -X_3\\
\star & -\frac{u_{lim}^2}{\omega_{max}^2} I & 0 \\
\star & \star & 0 \end{array} \right)\,+\,
\mathbf{He}(\mathcal{Q}^TG_3 \mathcal{P})<0}
\end{equation}
where ${\small \mathcal{P}=[I\quad  0 \quad -2\epsilon_3 I]}$ and
${\small \mathcal{Q}=[-I \quad -K^T \quad I]}$.
\end{Thm}
\medskip
\begin{proof}
Proof is similar to the Theorem 1. In this case $\mathcal{R}$ in
inequality (\ref{eq::Th1_2}) is
\begin{displaymath}{\small
\mathcal{R}=[2\bar{\epsilon_3} X_3 \quad 2\bar{\epsilon_3}
X_3K^T]}
\end{displaymath}
\end{proof}

\section{A multi-objective problem solved using the new dilated
MIs}\label{sec_multi-ob}

In this section, the new MIs are used to solve a multi-objective
problem. The problem considered here is that of controller design
for a system with bounded actuators exposed to a peak bound
disturbance, $w(t)^Tw(t)<w_{max}^2$. Nevertheless, the idea can be
applied to other multi-objective cases. The goal is to design a
controller which makes this system internally stable and
guarantees disturbance attenuation, while avoiding any violation
of saturation limit ($u_{lim}$). As mentioned in \cite{c3}, this
problem is a multi-objective problem, often yielding very
conservative results. Typically, the linear or low gain controller
for this problem is based on finding a controller (a state
feedback or dynamic output feedback compensator) which gives the
best $L_2$ performance while making sure that saturation limits
are not violated, by keeping the maximum of the control input
below the limit for all the points in the reachable set of the
closed loop system. We can formulate the proposed solution through
the following two algorithms: \smallskip
\begin{itemize}
\item \textbf{Algorithm 1} (Conventional approach): For common
$Q_1=Q_2=Q_3=Q>0$, minimize $\gamma_{con}$ in LMI
(\ref{eq::Q_Perf}) subject to MIs (\ref{eq::Q_Invariant}) and
(\ref{eq::Q_constraint}).

\smallskip

\item \textbf{Algorithm 2} (New approach using new dilated
MIs): For $X_1>0$, $X_2=X_3>0$ and common $G_1=G_2=G_3$, minimize
$\gamma_{new}$ in (\ref{eq::eq_perf}) subject to MIs
(\ref{eq::eq_Invariant}) and (\ref{eq::eq_constraint}).
\end{itemize}

\smallskip
In Algorithm 2, as we show later in both state feedback
and full order output feedback cases, common $G$ is needed to turn
the MIs into appropriate form to get a unique controller. However,
there is no obligation to use a common Lyapunov matrices for the
$L_2$ and invariant set inequalities. The reason to take $X_2=X_3$
in Algorithm 2 is that we are keeping the maximum of the
controller, below its limit, for all the points in the reachable
set of the closed-loop system. Here, the reachable set is
$x^TX_2^{-1}x<w_{max}^2$.

\smallskip
The following theorem states the advantage of the new
approach expressed through Algorithm 2 over the conventional one
obtained by Algorithm 1, by guaranteing better (or at least no
worse) $L_2$ gain estimate.

\medskip
\begin{Thm}[Multi-objective]\label{Thm_MO}
For the multi-objective saturation problem mentioned above,
Algorithm 2 with a common auxiliary variable, $G$, but with
non-common Lyapunov variables, always achieves an upper bound
estimate for the $L_2$ gain that is less than or equal to $L_2$
gain performance estimate achieved by Algorithm 1.
\end{Thm}
\medskip
\begin{proof}
If Algorithm 1 is solved, then Theorem \ref{TH_DMI_L2} implies
that there is a positive $\epsilon$ such that for any $\epsilon_1
< \epsilon$ by taking $X_1=G=G^T=Q$ and
$\gamma_{new}=\gamma_{con}$, we can satisfy (\ref{eq::eq_perf})
with the same closed-loop system derived by solving Algorithm 1.
Similarly, following the same argument, Theorems 2 and 3 guarantee
that there are $\epsilon_2$ and $\epsilon_3$ such that
$X_2=G=G^T=Q$ and $\alpha_{new}=\alpha_{con}$ satisfy
(\ref{eq::eq_Invariant}) and (\ref{eq::eq_constraint}). Therefore,
all the MIs in Algorithm 2 are feasible with
$\gamma_{new}=\gamma_{con}$ if we set $G=G^T=X_1=X_2=Q$ and use
the same closed-loop system, as obtained by Algorithm 1.

\smallskip
Therefore, any solution of Algorithm 1 can be achieved
by Algorithm 2, for a small enough $\epsilon$'s, without
exploiting the ability to use different Lyapunov matrices, which
could only improve the results.
\end{proof}

\smallskip
Recall that if MI (\ref{eq::eq_perf}) holds for some
$\epsilon_1$, it would hold for any $\epsilon<\epsilon_1$. Same
argument is true for MIs (\ref{eq::eq_Invariant}) and
(\ref{eq::eq_constraint}). Therefore, using the same $\epsilon$
still allows the results to be at least as good as those from
Algorithm 1. Therefore in Algorithm 2, to decrease the
computational cost of line search for $\epsilon_1$, $\epsilon_2$,
and $\epsilon_3$, we can use the same $\epsilon$ for all MIs. This
leads to some degree of conservatism. The conservatism is due to
the fact that the best result obtained by Algorithm 2 is not
necessarily preserved if we lower all or any of the $\epsilon_i$.
For best results, the $\epsilon$'s are allowed to vary
independently (see the numerical examples).

\smallskip
In Section \ref{sec_multi-ob_FSF},  Algorithms 1 and 2
are applied to state-feedback case and in Section
\ref{sec_multi-ob_OF}, these algorithms are used for full order
dynamic output feedback compensator. The result of these two
algorithms are compared through numerical examples.

\smallskip
\subsection{State Feedback Case}\label{sec_multi-ob_FSF}
In this section, state feedback controller synthesis is
considered. Therefore, $u=Kx_p$, and the closed-loop system in (3)
can easily be obtained with $x=x_p$.
For state feedback control, Algorithm 1 and 2 can be expressed in
the convex MI set as in lemmas below.

\smallskip
\begin{Lem}[\cite{c3}]\label{Lem_S1_FSF}
System (\ref{eq::CL_system}) with $w(t)^Tw(t) \le w_{max}^2$ and
state feedback K is internally stable, never saturates and has a
disturbance attenuation level $\gamma_{con}$ if there exist $Q>0$,
$Y$ and a positive constant $\alpha_{con}>0$ such that

\begin{displaymath}\label{eq::P1_1}{\small
\left( \begin{array}{ccc}
\Pi & B_1 & QC_1+Y^TD_{12}^T\\
\star & -\gamma_{con} I & D_{11}^T\\
\star & \star & -\gamma_{con} I
\end{array} \right)<0}
\end{displaymath}

\begin{displaymath}\label{eq::P1_2}{\small
\left( \begin{array}{cc}
\Pi+\alpha_{con} Q& B_1\\
\star & -\alpha_{con} I
\end{array} \right)<0}
\end{displaymath}

\begin{displaymath}\label{eq::P1_3}{\small
\left( \begin{array}{cc}
-Q& -Y^T\\
\star & -\frac{u_{lim}^2}{\omega_{max}^2} I
\end{array} \right)<0}
\end{displaymath}
where ${\small \Pi=AQ+QA^T+B_2Y+Y^TB_2}$. The variables in this
problem are $Q$, , $Y$, $\gamma_{con}$, and $\alpha_{con}$ where
$\alpha_{con}$ is searched through a line search. The controller
is given by $K=YQ^{-1}$.
\end{Lem}

\medskip
\begin{Lem}\label{Lem_S2_FSF}
System (\ref{eq::CL_system}) with $w(t)^T w(t) \le w_{max}^2$ and
state feedback K is internally stable, never saturates and has a
disturbance attenuation level $\gamma_{new}$ if there exist
$X_1>0$, $X_2>0$, $Y$, square matrix $G$, constant
$\alpha_{new}>0$, and small positive scalars $\epsilon_i<1$
$(i=1,2,3)$ such that

\setlength\arraycolsep{1pt}
{\small\begin{eqnarray}\label{eq::P2_1} \left( \begin{array}{cc}
X_1+\Pi+\Pi^T & B_1\\
\star & -\gamma_{new} I\\
\star &  \star \\
\star & \star
\end{array} \right.\qquad\qquad\qquad\qquad\qquad \nonumber \\
\left. \begin{array}{cc}
G^TC^T+Y^TD_{12}^T\quad & -X_1+G^T-2\epsilon_1 (\Pi) \\
D_{11}^T & 0 \\
-\gamma_{new} I & -2\epsilon_1(CG+D_{12}Y)\\
\star & -2\epsilon_1 (G^T+G)
\end{array} \right)<0\nonumber
\end{eqnarray}}
\smallskip
\setlength\arraycolsep{1pt}
{\small\begin{eqnarray}\label{eq::P2_2}
 \left( \begin{array}{cc}
X_2+\Pi+\Pi^T+\frac{\alpha_{new}}{2}(G+G^T)\qquad & B_1  \\
\star & -\alpha_{new} I  \\
\star &  \star
\end{array} \right. \nonumber \\
\left. \begin{array}{c}
-X_2+G^T-2\epsilon_2 (\Pi+\frac{\alpha_{new}}{2}G)\\
0 \\
-2\epsilon_2 (G^T+G)
\end{array} \right)<0\qquad\nonumber
\end{eqnarray}}
\smallskip
\setlength\arraycolsep{1pt}
{\small\begin{displaymath}\label{eq::P2_3} \left(
\begin{array}{cccc}
X_2-G-G^T & -Y^T & -X_2+G^T-2\epsilon_3 G\\
\star & -\frac{u_{lim}^2}{\omega_{max}^2} I & 2\epsilon_3 Y \\
\star & \star & -2\epsilon_3 (G^T+G) \end{array} \right)<0
\end{displaymath}}
where $\Pi=AG+B_2Y-\frac{1}{2}G$. The variables in this problem
are: Lyapunov matrices $X_1$ and $X_2$ as well as $\gamma_{new}$,
$Y$, $\alpha_{new}$ and $\epsilon_i\,\,(i=1,2,3)$ where
$\alpha_{new}$ and $\epsilon$'s are searched through line search.
The controller is given by $K=YG^{-1}$.
\end{Lem}

\smallskip
In this lemma, same slack variable G ($G_1=G_2=G_3=G$)
is used to make the MI set convex and to get a unique solution for
$K$.

\smallskip
\subsection{Numerical Example: State Feedback} In this section we use the results
of Section \ref{sec_multi-ob_FSF} in a numerical example. Consider
the example from \cite{c14} {\small\begin{equation}\label{example}
\setlength\arraycolsep{3pt} \left(\begin{array}{c|c|c} A & B_1 &
B_2\\\hline C_1 & D_{11} & D_{12} \end{array} \right)=\left(
\begin{array}{cccc|c|c}
0 & 0 & 1  &  0 & 0 & 0\\
0 & 0 & 0 &  1 & 0 & 0\\
-k & k &-f &f & 0 & 1\\
k & -k & f & -f & 1 & 0 \\\hline
0 & 1 &0 &0 & 0 & 0\\
0 & 0 &0 &0 & 0 & 0.01\end{array} \right)\end{equation}}with $k=2$
and $f=0.2$. The system is exposed to a peak bounded disturbance
with $w_{max}=5$. The controller limit is $u_{lim}=8$. Using
Algorithm 1, we get the following controller
$${\small K=\left[\begin{array}{cccc} -1.7970&   -0.7094&   -2.2916&   -2.1091\end{array}\right]}$$
which gives the minimum $\gamma_{con}^*=1.3038$. \noindent This
problem is also solved through Algorithm 2. To decrease the
computational cost, we used a same $\epsilon$ in all three dilated
MIs. Based on the simulations that we have done so far, the
variation of $\gamma_{new}$ with $\epsilon$ is bowl shaped.
Therefore, instead of a full line search on $\epsilon$, we used
the Golden-ratio method to identify the minimum performance,
$\gamma_{new}^*$, and the corresponding $\epsilon$. The result is
achieved very fast (9
iteration). 
The best performance that can be achieved by new method is
$\gamma_{new}^*=0.8345$ at $\epsilon=0.0802$, which is about
$36\%$
improvement, over the conventional method. 
 The controller associate with $\gamma_{new}^*$  is
$${\small K=\left[\begin{array}{cccc} -1.2732& -0.8923& -1.8967&  -1.8145
\end{array}\right]}$$



\medskip
Searching for three independent $\epsilon\,$'s can be a tedious
job, particularly if the number of objectives is rather large.
Here, to show the effect of independent $\epsilon\,$'s, we fix one
of the $\epsilon$'s, starting with the value obtained above (i.e.,
when all were set equal to one another),  and search for the other
two $\epsilon\,$'s, which are assumed to be equal. The best result
obtained is $\gamma^*=0.8104$ for $\epsilon_2=0.0802\,\,$ and
$\epsilon_1=\epsilon_3=0.1292$. Searching for three independent
$\epsilon$ leads on a slightly better result (about $3\%$), which
might be enhanced if a full -- or more thorough -- search is done,
but the computation burden will be significant.

\medskip


\subsection{Output Feedback Case} \label{sec_multi-ob_OF} In this section, we
use the dilated MIs in designing a full order dynamic output
feedback compensator for the same saturation problem. Therefore,
the controller is
\begin{equation}  \label{eq::output_feedback}{\small
\left \{\begin{array}{l}
\dot{x_c}  = A_cx_c+B_cy \\
u  =  C_cx_c \end{array} \right.}
\end{equation}
where $A_c$ is of the same order as the system matrix $A$. To
reduce clutter of the equations, we dropped $D_c$. Applying this
controller, the closed-loop system in (3) can easily be obtained
with $x =[x_p^T\,\,\, x_c^T]^T$.

\medskip
In output feedback case, additional complications arise since a
variety of transformations and manipulations are needed to set the
problem into a convex search, often requiring auxiliary variables
instead of the compensator matrices ($A_c$, $B_c$ and $C_c$). The
approach is reasonably well known and can be found in a variety of
references (e.g., \cite{c5} and \cite{c15} among many).  Here,
most of the technical details are omitted but can be found in the
references mentioned, though the outline is based on the approach
used in \cite{c5} and \cite{c15}.  Here, for the  Lyapunov matrix
in Algorithm 1, we use the structure
$${\small P=\left( \begin{array}{cc}
Y & -Y\\
-Y & S^{-1}+Y
\end{array} \right)>0}$$
which, as discussed in \cite{c15},  can be done without any loss
of generality. Therefore, Algorithm 1 can be expressed through the
following lemma.

\medskip
\begin{Lem}[\cite{c3}]\label{Lem_S1_OF}
System (\ref{eq::CL_system}) with disturbance $w(t)$ satisfying
$w(t)^T w(t) \le w_{max}^2$ and with a output feedback compensator
(\ref{eq::output_feedback}) is internally stable, never saturates
and has a disturbance attenuation level $\gamma_{con}$ if there
exist $Y>0$, $X>0$, and general matrices $L$, $F$ and $E$ and a
constant $\alpha_{con}>0$ such that

\begin{displaymath}\label{eq::OP_P1_1}{\small
\left( \begin{array}{cccc}
\Pi & A+L^T & B_1 & X^TC_1^T+F^TD_{12}^T\\
\star & \Lambda & YB_1+ED_{21}&C1\\
\star & \star & -\gamma_{con} I & D_{11}\\
\star & \star & \star & -\gamma_{con} I
\end{array} \right)<0}
\end{displaymath}
\begin{displaymath}\label{eq::OP_P1_2}{\small
\left( \begin{array}{ccc}
\Pi+\alpha_1 X& A+L^T+\alpha_{con} I & B_1\\
\star & \Lambda+\alpha_{con} Y & YB_1+ED_{21}\\
\star &  \star & -\alpha_{con} I
\end{array} \right)<0}
\end{displaymath}
\begin{displaymath}\label{eq::OP_P1_3}{\small
\left( \begin{array}{ccc}
-X& -I & -F^T\\
\star & -Y & 0\\
\star & \star & -\frac{u_{lim}^2}{\omega_{max}^2} I
\end{array} \right)<0}
\end{displaymath}
where ${\small\Pi=\mathbf{He}(AX+B_2Y)}$ and
${\small\Lambda=\mathbf{He}(A^TY+EC_2)}$. Then for ${\small
S=X-Y^{-1}}$, one representation of the controller matrices is:
$${\small C_c=FS^{-1},\quad
B_c=-Y^{-1}E}$$
$${\small A_c=(A-B_cC_2)XS^{-1}+B_2C_c-Y^{-1}LS^{-1}}$$
\end{Lem}

\bigskip
For dilated MIs, also we need to do some manipulations to
expand them into appropriate convex or near convex forms. We start
with dilated MI for $L_2$ gain (\ref{eq::eq_perf}). As in the
state feedback multi-objective solution, eventually we are going
to use common $G$, the key slack variable introduced by dilation.
As a result, for simplification, the index $i$ of
$G_i\,\,(i=1,2,3)$ is dropped. Let us call $G^{-1}=H$. Note that
based on the structure of dilated MIs, $G$ is invertible. By pre-
and post-multiplying (\ref{eq::eq_perf}) by
\begin{displaymath}{\small
Diag \left [ \begin{array}{cccc} T^TH^T & I & I & T^T H^T
\end{array} \right]}
\end{displaymath}
and its transpose, respectively, we obtain
\begin{equation}\label{eq::OF_perf}{\small
 \left( \begin{array}{cccc}
\Phi_{11} & T^TH^TB_{cl}  &T^TC_{cl}^T\quad & \Phi_{14} \\
\star & -\gamma_2 I  & D_{cl}^T & 0 \\
\star &  \star & -\gamma_2 I & -2\epsilon_1(C_{cl}T)\\
\star & \star & \star & -2\epsilon_1 \Phi_{44}
\end{array} \right)<0}
\end{equation}
where
$${\small \left. \begin{array}{c} \Phi_{11}=T^TH^TX_1HT+\mathbf{He}(T^TH^T
A_{cl}T)-\frac{1}{2}\Phi_{44}\\
 \Phi_{14}=-T^T(H^TX_1H +H-2\epsilon_1
(H^TA_{cl}-\frac{1}{2}H^T))T\\
\Phi_{44}=\mathbf{He}(T^TH^TT)
\end{array}\right. }$$
Here, $T$ is an auxiliary matrix used for the additional
transformations that are needed in the output feedback synthesis
problem.

\medskip
We partition $H$ and $G$ into the following forms with each
submatrix having the dimension $n\times n$
\begin{equation}\label{H}{\small
H=\left(\begin{array}{cc}
H_{11} & H_{12}\\
H_{21} & H_{22}\end{array}\right), \quad G=\left(\begin{array}{cc}
G_{11} & G_{12}\\
G_{21} & G_{22}\end{array}\right).}
\end{equation}
As mentioned before, $G$ is invertible, so the sub-matrices of
this matrix are also invertible (by invoking a small perturbation
if necessary \cite{c16}). Now, to turn dilated MIs into near
convex form, let's consider the the auxiliary transformation
matrix $T$ as
$${\small T=\left(\begin{array}{cc}
G_{11} & I\\
G_{21} & 0\end{array}\right).}
$$
To simplify notations, we call $G_{11}\equiv R$ and $H_{11}\equiv
Y$. Finally, let us call
$${\small T^THX_1HT=M_1= \left(\begin{array}{cc}
M_{11} & M_{12}\\
M_{12}^T & M_{13}\end{array}\right);\, (M_1=M_1^T)}$$ and
$${\small F=C_cG_{21},\quad E=H_{21}^TB_c,\quad V=R^TY+G_{21}^TH_{21}}$$
$${\small L=(Y^TA+EC_2)R+Y^TB_2F+H_{21}^TA_cG_{21}}$$ These are the new
variables introduced to turn this MIs into convex form (similar to
$G$, $F$, $L$ used in the standard approach. Also to save space,
for $i=1,2$, we use
$${\small \Pi=AR+B_2F-\frac{1}{2}R,\quad
\Lambda=YA+EC_2-\frac{1}{2}Y^T}$$
$${\small\Omega=-\frac{1}{2}(I+V)+A+L^T,~
\Delta_i=-M_{i1}+R^T-2\epsilon_i\Pi}$$
$${\small\Sigma_i=-M_{i3}+Y-2\epsilon_i\Lambda,~
\Gamma_i=-M_{i2}+V-2\epsilon_i(A-\frac{1}{2}I)}$$
$${\small\Upsilon_i=-M_{i2}^T+I-2\epsilon_i(L-\frac{1}{2}V^T)}$$

Naturally, the same manipulation can be conducted on MI
(\ref{eq::eq_Invariant}) and (\ref{eq::eq_constraint}), albeit
with a bit less clutter. Considering above matrices and
definitions, the entries of the three MIs can be expanded in
detail as in the following lemma which expresses Algorithm 2 for
dynamic output feedback compensator.

\medskip
\begin{Lem}\label{Lem_S2_OF}System (\ref{eq::CL_system})
with disturbance $w(t)$ satisfying $w^T(t) w(t) \le w_{max}^2$ and
with a output feedback compensator (\ref{eq::output_feedback}) is
internally stable, never saturates and has a disturbance
attenuation level $\gamma_{new}$ if there exist square matrices
$R$, $S$, $V$ and symmetric matrices $Y$, $M_1$, $M_2$
\setlength\arraycolsep{3pt}
$${\small M_1= \left(\begin{array}{cc}
M_{11} & M_{12}\\
M_{12}^T & M_{13}\end{array}\right),\quad M_2=
\left(\begin{array}{cc}
M_{21} & M_{22}\\
M_{22}^T & M_{23}\end{array}\right)}$$ and general matrices
$L$,$E$, $F$, and constant $\alpha_{new}>0$, and small positive
value $\epsilon_i<1$ $(i=1,2,3)$ such that
{\setlength\arraycolsep{0.1pt}}
$${\small\left(
\begin{array}{ccc}
M_{11}+\mathbf{He}(\Pi) & M_{12}+\Omega  & B_1 \\
\star & \,M_{13}+\mathbf{He}(\Lambda)\, & Y^TB_1+ED_{21} \\
\star &  \star & -\gamma_{new} I \\
\star & \star &\star  \\
\star & \star &\star \\
\star & \star &\star
\end{array} \right.}
$${\setlength\arraycolsep{0.1pt}}
$${\small \left. \begin{array}{ccc}
R^TC_1^T+F^TD_{12}^T & \Delta_1 & \Gamma_1 \\
C_1^T & \Upsilon_1 & \Sigma_1\\
D_{11}^T & 0 & 0 \\
-\gamma_{new} I&  -2\epsilon_1(C_1R+D_{12}F) & -2\epsilon_1C_1\\
\star&  -2\epsilon_1 \mathbf{He}(R) & -2\epsilon_1 (I+V)\\
\star &  \star & -2\epsilon_1 \mathbf{He}(Y)
\end{array} \right)<0\quad}$$
\medskip {\setlength\arraycolsep{0.1pt}}
$${\small\left(
\begin{array}{cc}
M_{21}+\mathbf{He}(\frac{\alpha_{new}}{2}R+\Pi) & M_{22}+\frac{\alpha_{new}}{2} (I+V)+\Omega \\
\star & M_{23}+\mathbf{He}(\frac{\alpha_{new}}{2}Y+\Lambda)\\
\star &  \star \\
\star &\star\\
\star &\star\\
\end{array} \right.}$${\setlength\arraycolsep{0.005pt}}
$${\small\left. \begin{array}{ccc}
B_1 & \Delta_2-2\epsilon_2(\frac{\alpha_{new}}{2}R)\quad & \Gamma_2-2\epsilon_2(\frac{\alpha_{new}}{2}I) \\
Y^TB_1+ED_{21} & \Upsilon_2-2\epsilon_2(\frac{\alpha_{new}}{2}V^T)
&
\Sigma_2-2\epsilon_2(\frac{\alpha_{new}}{2}Y^T)\\
-\alpha_{new} I & 0 & 0 \\
\star & \mathbf{He}(-2\epsilon_2 R) & -2\epsilon_2 (I+V)\\
\star & \star & \mathbf{He}(-2\epsilon_2 Y)
\end{array} \right)<0\,\,}$$
\smallskip
$${\small\left(
\begin{array}{ccc}
M_{21}-\mathbf{He}(R) & M_{22}-I-V & -F^T \\
\star & M_{23}-\mathbf{He}(Y) & 0 \\
\star &  \star & -\frac{u^2_{max}}{\omega^2_{max}}I \\
\star &\star & \star \\
\star &\star &\star
\end{array} \right.}$$
$${\small\left. \begin{array}{cc} -M_{21}+R^T+2\epsilon_3R\quad
&-M_{22}+V+2\epsilon_3I\\
-M_{22}^T+I+2\epsilon_3 V^T & -M_{23}+Y+2\epsilon_3Y^T \\
F & 0\\
\mathbf{He}(-2\epsilon_3 R) & -2\epsilon_3 (V+I)\\
\star & \mathbf{He}(-2\epsilon_3 Y)\\
\end{array} \right)<0}$$

\normalsize
\medskip
Then for invertible $G_{21}$ and $H_{21}$ deduced by
matrix factorization of $G_{21}^TH_{21}=V-R^TY$ as mentioned in
\cite{c17}, the controller matrices are as follows:
$${\small C_c=FG_{21}^{-1},\quad B_c=H_{21}^{-T}E}$$
$${\small A_c=H_{21}^{-T}(L-Y^TAR+EC_2R+YB_2F)G_{21}^{-1}}$$
\end{Lem}

\medskip
As mentioned above, invertible $G_{21}$ and $H_{21}$ can
be deduced by matrix factorization of $G_{21}^TH_{21}=V-R^TY$. As
indicated in \cite{c17}, this deduction is always possible and if
necessary we can use perturbation. In our numerical examples, when
possible, we obtain $G_{21}$ and $H_{21}$ by picking $H_{21}= Y$
and therefore having $G_{21}=Y^{-T}V^T-R$.

\smallskip
As before, to avoid excessive computational cost in the
search for $\epsilon_1$, $\epsilon_2$, and $\epsilon_3$ in Lemma
(\ref{Lem_S2_OF}), we can use same $\epsilon$ for all MIs.

\medskip
\subsection{Numerical example: Output Feedback} Consider the same numerical
example as for the state feedback, this time with $k=0.4$,
$f=0.04$, $u_{lim}=100$ and with $y^T=[x_1\quad x_2]$.

\smallskip
This system is expected to withstand disturbances with
peak bound
of $w_{max}=5$. 
Using the conventional approach, the resulting controller has
$A_c$ and $B_c$ matrices of order $10^{5}$ and
$${\small C_c|_{con}=\left[
\begin{array}{cccc} -12.5082 & -18.5711 & -5.8239 & -42.4745
\end{array} \right]}$$
The induced $L_2$ gain of this system under this controller is
$\gamma^*_{con}=1.7636$.

\smallskip
This problem is also solved with the new dilated matrix
inequalities. We start with the same $\epsilon$ in all three
dilated matrix inequalities, and we use Golden-ratio method in
search over $\epsilon$. This, after 10 iteration, leads to
$\gamma_{new}^*=1.5029$ at $\epsilon=0.0231$, which is about
$15\%$ improvement, over the conventional method. The
corresponding controller has $A_c$ and $B_c$ matrices of order
$10^4$ and
$${\small C_c|_{new}=\left[ \begin{array}{cccc} -11.5857 &
-17.0741 & -5.5596 & -38.7530
\end{array} \right]}$$

\smallskip
Next, we try independently varying $\epsilon$'s. To
limit the cost of searching for three independent $\epsilon$'s, we
have fixed one of the $\epsilon$'s in a certain value and search
for the other two $\epsilon$'s which are assumed to be equal. The
best result obtained is $\gamma^*=1.2746$ for $\epsilon_2=0.0181$
and $\epsilon_1=\epsilon_3=0.1618$. Thus, letting $\epsilon$'s
vary independently has a significant effect at the results, an
improvement of about $27\%$ over the conventional method of
solving the problem ($\gamma^*_{con}=1.7636$).


\medskip
\section{CONCLUSIONS}

We presented new dilated matrix inequalities for Bounded Real MI,
invariant set MI and constraint MI. The structure of these dilated
MIs, in which system matrices are separated from Lyapunov
matrices, allows us to use different Lyapunov matrices for
different objective in multi-objective problems or different
parameter values in robust synthesis problems. The new approach is
guaranteed to achieve results which are better or equal to the
ones obtained from the standard multi-objective setting. The
synthesis results, for both state feedback and output feedback
problems, are demonstrated through an example.


\end{document}